\documentclass[a4paper, oneside, 10pt]{article}
\usepackage[english]{babel}
\usepackage{amsfonts}
\usepackage{amsmath}
\usepackage{amssymb}
\usepackage{graphics}
\usepackage{amsrefs}
\usepackage{latexsym}
\begin{document}

\title{{\bf{\Large{On some Bayesian nonparametric estimators for species richness under two-parameter Poisson-Dirichlet priors}}}\footnote{{\it AMS (2000) subject classification}. Primary: 60G58. Secondary: 60G09.}}
\author{\textsc {Annalisa Cerquetti}\footnote{Corresponding author. Annalisa Cerquetti, Dipartimento di Studi Geoeconomici, Sapienza Universit\`a di Roma, Via del Castro Laurenziano, 9, 00161 Roma.
E-mail: {\tt annalisa.cerquetti@uniroma1.it}}\\
  \it{\small Department of Geoeconomic, Linguistic, Statistical and Historical Studies}\\ \it {\small SAPIENZA University of Rome, Italy }}
\newtheorem{teo}{Theorem}
\date{\today}
\maketitle{}

\begin{abstract}

We present an alternative approach to the Bayesian nonparametric analysis of  conditional species richness under two-parameter Poisson Dirichlet priors. We rely on a known characterization by {\it deletion of classes} property and on results for Beta-Binomial distributions. Besides leading to simplified and much more direct proofs, our proposal provides a new scale mixture representation of the conditional asymptotic law.
\end{abstract}

\section{Introduction}

In Favaro et al. (2009) explicit expressions for Bayesian nonparametric estimators for conditional species richness under two-parameter Poisson-Dirichlet priors have been derived to deal with the problem of prediction when the size of the additional sample tends to be very large. The paper also investigates the asymptotic behavior of this quantity in order to obtain asymptotic highest posterior density intervals for the estimates of interest.  Despite referring to the Bayesian nonparametric treatment of {\it conditional Gibbs structures} as introduced in Lijoi et al. (2007, 2008), the proofs are somehow cumbersome and do not resort to previously established properties of these structures nor to some specific available results for the two-parameter Poisson-Dirichlet family. \\

Here we show how the results in Favaro et al. (2009) may be derived by a much more direct and simpler approach, resorting to the {\it deletion of classes} property of the two-parameter Poisson-Dirichlet model (Pitman, 2003) and to known properties of the Beta-Binomial distribution. Moreover, as a by product, we obtain a new scale mixture representation for the limit law of the conditional species richness which differs from that derived in Favaro et al. (2009). \\

Notice that, to make the paper easily readable to those unfamiliar with the Bayesian treatment of exchangeable Gibbs partitions, we adopt Pitman's (2006) notation.  A preliminary rephrasing of the Lijoi et al. (2007, 2008) approach in terms of Pitman's theory  may be found in Cerquetti (2008), where even the relationship between conditional Gibbs structures and the operation of {\it deletion of classes} has been first pointed out. For the sake of clarity and to make the paper self-contained, we open each section with known results of which we will make use throughout the paper.

\section{Some preliminaries on the two-parameter Poisson-Dirichlet partition model}

The two-parameter $(\alpha, \theta)$  Poisson-Dirichlet distribution for $\alpha \in [0,1)$ and $\theta \geq -\alpha$ is a model  for random partitions (Pitman and Yor, 1997), which belongs to the class of exchangeable Gibbs partitions of type $\alpha \in [-\infty, 1)$ as defined in Gnedin and Pitman (2006). This class is characterized by an exchangeable partition probability function of the form
$$
p(n_1,\dots, n_k)=V_{n,k} \prod_{j=1}^k (1-\alpha)_{n_j-1},
$$
where $(a)_b=a(a+1)\cdots(a+b-1)$ are rising factorials, with weights $V_{n,k}$ satisfying the backward recursion $V_{n,k}=(n -k\alpha)V_{n+1,k} +V_{n+1, k+1}$.
The $(\alpha, \theta)$ Poisson-Dirichlet distribution is well-known to arise for 
\begin{equation}
\label{pesi}
V_{n,k}=\frac{(\theta +\alpha)_{k-1\uparrow \alpha}}{(\theta +1)_{n-1}}, 
\end{equation}
where $(x)_{s\uparrow \alpha}$ stands for generalized rising factorials $(x)_{s \uparrow \alpha}= (x)(x +\alpha)(x+2\alpha)\cdots (x +(s-1)\alpha)$.
This model has been largely studied in the last twenty years (see e.g. Perman et al. 1992, Pitman, 1995, 1996a, 1996b, Pitman, 2003) and a lot of results are available for it. Here we just recall few of them that we are going to exploit in the following.  A general reference is Pitman (2006). \\

For $S_{n,k}^{-1, -\alpha}$ the generalized Stirling numbers of the first kind (see e.g. Hsu \& Shiue, 1998), the law of the number of blocks $K_n$ observed in an $n$-sample for $k \in \{1, \dots, n\}$ is given by
\begin{equation}
\label{blocksPD}
\mathbb{P}_{\alpha,\theta}(K_n=k)=\frac{(\theta +\alpha)_{k-1 \uparrow \alpha}}{(\theta +1)_{n-1}} S_{n,k}^{-1, -\alpha},
\end{equation}
 with expected value equal to
\begin{equation}
\label{primoPD}
\mathbb{E}_{\alpha, \theta}(K_n)=\frac{(\theta +\alpha)_{n}}{\alpha(\theta +1)_{n-1}} - \frac{\theta}{\alpha}.
\end{equation}\\
A general expression for the moments of any order of $K_n$ has been obtained in Yamato and Sibuya (2000) and Pitman (1996b) in terms of {\it non-central} Stirling numbers of the second kind $S_{r,j}^{0,1,\theta/\alpha}$, 
and corresponds to
\begin{equation}
\label{momer}
\mathbb{E}_{\alpha, \theta}(K_n^{r})= \sum_{j=0}^{r} (-1)^{r-j}\left(\theta/\alpha +1 \right)_{j} S_{r,j}^{0, 1, \theta/\alpha} \frac{(\theta +j\alpha +1)_{n -1}}{(\theta +1)_{n -1}}.
\end{equation}
\\
Now, {\it conditional Gibbs structures} have been introduced as tools for a Bayesian nonparametric approach to species sampling problems in Lijoi et al. (2007, 2008). 
In this setting, given a sample $(X_1, \dots, X_n)$,  with $(n_1, \dots, n_k)$ the vector of observed multiplicities of each species represented,  interest typically lies in the law of the number $K_m$ of different species observed in an additional $m$-sample $(X_{n+1}, \dots, X_{n+m})$. The general form of this distribution 
have been first obtained in Lijoi et. al. (2007, cfr. Proposition 1.) by combinatorial arguments, and may be expressed in terms of {\it non-central} generalized Stirling numbers of the first kind (cfr. Cerquetti, 2008, cfr. Eq. 32) as follows
\begin{equation}
\label{newblocks}
\mathbb{P}(K_m=k^*|n_1, \dots, n_k)=\mathbb{P}(K_m=k^*|K_n=k)=\frac{V_{n+m, k+k^*}}{V_{n,k}} S_{s,k^*}^{-1, -\alpha, -(n -k\alpha)},
\end{equation}
for $k^* \in \{0, 1, \dots, m\}$. This formula specializes under the $(\alpha, \theta)$ Poisson-Dirichlet model as
\begin{equation}
\label{blockskm}
\mathbb{P}_{\alpha, \theta}(K_m=k^*|K_n=k)=\frac{(\theta +k\alpha)_{k^* \uparrow \alpha}}{(\theta +n)_{m}} S_{m, k^*}^{-1, -\alpha, -(n -k\alpha)}, 
\end{equation}
whose expected value, given by
\begin{equation}
\label{att}
\mathbb{E}_{\alpha, \theta}(K_m|K_n=k)=\sum_{k^*=0}^{m} k^* \frac{(\theta +k\alpha)_{k^* \uparrow \alpha}}{(\theta +n)_{m }} S_{m, k^*}^{-1, -\alpha, -(n -k\alpha)},
\end{equation}
plays the role of a Bayesian estimator for $K_m$. 

\section{Conditional analysis for species richness under two-parameter Poisson-Dirichlet priors}

Favaro et al. (2009) move from the need of an alternative expression for (\ref{att}) to reduce the computational effort needed to calculate both (\ref{att}) and Bayesian estimators for related quantities of interest in species sampling problems. These basically sum up to the {\it discovery probability}, the probability to discover a new species at the $(n+m+1)$th draw without observing the $m$ intermediate records,  and the {\it sample coverage},  the proportion of species represented in a sample of given size featuring a certain number of distinct species.\\

Here is our approach to the problem. Let $S_m$ be the number of observations in the additional $m$-sample belonging to new species, with values in $\{0,\dots,m\}$. By the basic rules of conditional probability we can always write (\ref{newblocks}) as
$$
\mathbb{P}(K_m = k^* | K_n =k)= \sum_{s=0}^{m} \mathbb{P}(K_m=k^*, S_m=s|K_n=k)=$$
\begin{equation}
\label{mix}
= \sum_{s=0}^{m} \mathbb{P}(K_m=k^* |K_n=k, S_m=s) \mathbb{P} (S_m=s |K_n=k).
\end{equation}\\
The general form of $\mathbb{P} (S_m=s |K_n=k)$ for Gibbs partitions of type $\alpha \in [0,1)$  has been derived in Lijoi et al. (2008, cfr. Eq. (11)),  and expressed in terms of generalized Stirling numbers  as in Cerquetti (2008) is given by 
\begin{equation}
\label{obsinnew3}
\mathbb{P}(S_m=s|K_n=k)= \frac{1}{V_{n,k}}{m \choose s}(n-k\alpha)_{m-s}\sum_{k^*=0}^s {V_{n+m, k+k^*}} S_{s,k^*}^{-1,-\alpha}.
\end{equation}
This formula specializes under the $(\alpha, \theta)$-Poisson-Dirichlet model as follows. First notice that
$$
\frac{V_{n+m,k}}{V_{n,k}}=\frac{1}{(\theta +n)_{m}},
$$
then, by means of the multiplicative property of generalized rising factorials and the definition of generalized Stirling numbers as connection coefficients, the sum in (\ref{obsinnew3}) reduces to
$$
\frac{1}{(\theta +1)_{n+m-1}}\sum_{k^*=0}^s (\theta +\alpha)_{k+k^*-1\uparrow \alpha} S_{s,k^*}^{-1, -\alpha}=
\frac{(\theta +\alpha)_{k-1 \uparrow \alpha}}{(\theta +1)_{n+m-1}}\sum_{k^*=0}^s  (\theta + k \alpha)_{k^* \uparrow \alpha} S_{s, k^*}^{-1, -\alpha}=
$$
$$= \frac{(\theta +\alpha)_{k-1 \uparrow \alpha}}{(\theta +1)_{n+m-1}} (\theta +k\alpha)_{s}.
$$
It follows that (\ref{obsinnew3})  specializes  under the $(\alpha, \theta)$ model as\\\
$$
\mathbb{P}_{\alpha, \theta}(S_m=s|K_n=k)
={m \choose s}(n -k \alpha)_{m-s} \frac{V_{n+m,k}}{V_{n,k}}(\theta +k\alpha)_s=$$
\begin{equation}
\label{polya} 
={m \choose s}\frac{(n-k\alpha)_{m-s}(\theta +k\alpha)_{s}}{(\theta +n)_{m}},
\end{equation}
which is a {\it Beta-Binomial  $(m, \theta +k\alpha, n- k\alpha)$ distribution} with expected value
$$
\mathbb{E}_{\alpha, \theta}(S_m|K_n=k)=m\frac{(\theta +k\alpha)}{\theta +n}.\\\\
$$
\\
{\bf Remark 1.} In Lijoi et al. (2008) an analogous derivation of (\ref{polya}) in terms of {\it generalized factorial coefficients} (see Charalambides, 2005) is in Example 3.2. Nevertheless the relationship with Beta-Binomial distributions is not highlighted, (e.g. in deriving the expected value they resort to a general formula in Proposition 2.).  Notice that Beta-Binomial distributions (see e.g. Johnson and Kotz, 1977, 2005) can be seen as a generalization to non-integer parameters of {\it P\'olya urn distributions} of parameters $(a, b, c=1)$, for $a$ and $b$ the initial composition of the urn and $c$ the number of balls of the same color replaced in the urn with the ball observed. Asymptotic results for P\'olya distributions extend to Beta-Binomial models, something that we will exploit in the following sections. \\

As for $\mathbb{P}(K_m=k^* |K_n=k, S_m=s)$ this is the law of the number of blocks for a conditional Gibbs structure as defined in Lijoi et al. (2008, Prop. 3). As shown in Cerquetti (2008, Section 4.1),  the operation of conditioning to the number $s$ of observations in the new blocks is equivalent to conditioning to the number $m-s$ of observations in old blocks, i.e. to the vector $(m_1, \dots, m_k)$ and corresponds to the operation of {\it deletion of the first k classes} as defined in Pitman (2003). \\\\
{\bf Definition 2} [Deletion of classes, Pitman (2003)] Given a random partition $\Pi$ of $\mathbb{N}$, the operator {\it deletion of the first $k$ classes} is as follows: First let $\Pi_k^*$ be the restriction of $\Pi$ to $H_k:=\mathbb{N}-G_1-\cdots-G_k$ where $G_1,\dots, G_k$ are the first $k$ classes of $\Pi$ in order of their least elements, then derive $\Pi_k$ on $\mathbb{N}$ from $\Pi^*_k$ on $H_k$ by renumbering the points of $H_k$ in increasing order. \\\\
Pitman (2003) shows  that this operation characterizes  the two-parameter Poisson Dirichlet family of distributions for $\alpha \in (0,1)$ in that produces a Gibbs partition still belonging to the Poisson-Dirichlet class with updated parameter $(\alpha, \theta +k\alpha)$ (see Gnedin et al. 2009 for recent results and a comprehensive treatment of the topic). It follow that, by (\ref{blocksPD}) 
\begin{equation}
\label{deletion}
\mathbb{P}_{\alpha, \theta}(K_m=k^* |K_n=k, S_m=s)=\mathbb{P}_{\alpha, \theta +k\alpha}(K_s=k^*)=\frac{(\theta +k\alpha +\alpha)_{k^*- 1 \uparrow \alpha}}{(\theta +k\alpha +1)_{s-1}} S_{s, k^*}^{-1, \alpha}.
\end{equation}

\bigskip 

\noindent Notice that by (\ref{mix}), (\ref{polya}) and (\ref{deletion})  it is possible to reobtain (\ref{blockskm}) as follows
$$
\mathbb{P}_{\alpha, \theta}(K_m = k^* | K_n =k)= \frac{(\theta +k\alpha +\alpha)_{k^*- 1 \uparrow \alpha}}{(\theta +n)_{m}}\sum_{s=k^*}^m {m \choose s} \frac{(n-k\alpha)_{m-s}(\theta +k\alpha)_{s}}{(\theta +k\alpha +1)_{s-1}} S_{s, k^*}^{-1, \alpha}=
$$
by the definition of generalized rising factorial $(x)_{n \uparrow h}:=h^n(x/h)_{n \uparrow 1}$ reduces to
$$
=\frac{(\theta +k\alpha)}{(\theta +n)_m}\frac{\alpha^{k^*-1}\Gamma(\theta/\alpha +k +k^*)}{\Gamma(\theta/\alpha +k +1)}\sum_{s=k^*}^{m} {m \choose s}(n- k\alpha)_{m-s} S_{s,k^*}^{-1, -\alpha}=
$$
and by the definition of non-central generalized Stirling number as connection coefficients yields
$$
=\frac{(\theta +k\alpha)_{k^* \uparrow \alpha}}{(\theta +n)_m} S_{s,k^*}^{-1, -\alpha, -(n-k\alpha)}.
$$
Now we show how the approach described in the present section applies to the study of expected value, moments and the asymptotic behaviour  of $(K_m|K_n=k)$ under $(\alpha, \theta)$ Poisson-Dirichlet model.

\subsection{Moments}
By the mixture representation (\ref{mix}), and the simplification induced by the deletion of classes property, the moments of any order for the number of species in the additional sample conditional on the basic sample are given by:
\begin{equation}
\label{mom}
\mathbb {E}_{\alpha, \theta }(K_m^{r}|K_n=k)=\sum_{s=0}^{m} \mathbb{E}_{\alpha, \theta}(K_m^{r}|S_m=s, K_n=k) \mathbb{P}_{\alpha, \theta}(S_m=s|K_n=k).
\end{equation}
For $r=1$ deriving $\mathbb{E}_{\alpha, \theta}(K_m^{r}|S_m=s, K_n=k)$ is just a matter of specializing (\ref{primoPD}),
\begin{equation}
\label{mprimoks}
\mathbb{E}_{\alpha, \theta}(K_m|S_m=s, K_n=k)=\mathbb{E}_{\alpha, \theta+k\alpha}(K_s)=\frac{(\theta +k\alpha +\alpha)_s}{\alpha(\theta +k\alpha +1)_{s-1}}-\frac{\theta +k\alpha}{\alpha}.
\end{equation}
For $r >1$ specializing (\ref{momer}) yields 
$$
\mathbb{E}_{\alpha, \theta}(K_m^r|S_m=s, K_n=k)=
$$
\begin{equation}
\label{rmag}
\mathbb{E}_{\alpha, \theta+k\alpha}(K_s^{r})= \sum_{j=0}^{r} (-1)^{r-j}\left(\frac {\theta +k\alpha +\alpha}{\alpha} \right)_{j} S_{r,j}^{0,1, {(\theta +k\alpha})/{\alpha}} \frac{(\theta +k\alpha +j\alpha +1)_{s -1}}{(\theta +k\alpha +1)_{s -1}}.
\end{equation}\\\\
We are now in a position to prove Favaro et al. (2009) Proposition 1. in a much more direct and simple fashion.\\\\
{\bf Proposition 3.} Under the $(\alpha, \theta)$ Poisson-Dirichlet model an explicit expression for the expected value of $K_m$ conditioned to the number of blocks $K_n$  observed in the basic $n$-sample is as follows
\begin{equation}
\label{mioprimo}
\mathbb{E}_{\alpha, \theta}(K_m|K_n=k)= \left(\frac{\theta +k\alpha}{\alpha}\right) \left[\frac{(\theta +\alpha +n)_m}{(\theta +n)_m}-1\right].
\end{equation}
{\it Proof}: By (\ref{polya}), (\ref{mom}) and (\ref{mprimoks}) 
$$
 \mathbb{E}_{\alpha, \theta}(K_m|K_n=k)=\sum_{s=0}^{m} \mathbb{E}_{\alpha, \theta +k\alpha}(K_s) \mathbb{P}_{\alpha, \theta}(S_m=s| K_n=k)=
$$
$$
=\sum_{s=0}^{m} {m \choose s} \left[\frac{(\theta +k\alpha +\alpha)_s}{\alpha(\theta +k\alpha +1)_{s-1}}-\frac{\theta +k\alpha}{\alpha}\right] \frac{(n-k\alpha)_{m-s}(\theta +k\alpha)_s}{(\theta +n)_m}=
$$
$$
=\frac{1}{(\theta +n)_m} \sum_{s=0}^{m} {m \choose s} \left[ \frac{\theta +k\alpha}{\alpha}\left( \frac{(\theta +k\alpha +\alpha)_s}{(\theta +k\alpha)_s} -1\right)\right](\theta +k\alpha)_s (n -k\alpha)_{m-s}=
$$
$$=\left(\frac{\theta +k\alpha}{\alpha}\right)\frac{1}{(\theta +n)_m} \left[\sum_{s=0}^{m} {m \choose s} (\theta +k\alpha +\alpha)_s
(n-k\alpha)_{m-s} - (\theta +n)_m\right]=
$$
$$
=\left(\frac{\theta +k\alpha}{\alpha}\right)\left[ \frac{(\theta +\alpha +n)_m}{(\theta +n)_m} -1\right]. 
$$\hspace{14.7cm} $\square$\\\\
{\bf Proposition 4.} Under the $(\alpha, \theta)$ Poisson-Dirichlet model an explicit expression for the moments of any order for $(K_m|K_n=k)$ is given by
\begin{equation}
\mathbb{E}_{\alpha, \theta}(K_m^{r}|K_n=k)= \sum_{j=0}^r (-1)^{r-j} \left( \frac{\theta +k\alpha}{\alpha}\right)_j S_{r,j}^{0, 1, (\theta+k\alpha)/\alpha} \frac{(\theta +n +j\alpha)_m}{(\theta +n)_m}.
\end{equation}
{\it Proof}: By (\ref{polya}), (\ref{mom}) and (\ref{rmag})  
$$
\label{successivi}
\mathbb{E}_{\alpha, \theta}(K_m^{r}|K_n=k)=
$$
$$
= \sum_{s=0}^m {m \choose s} \sum_{j=0}^{r} (-1)^{r-j}\left(\frac {\theta +k\alpha +\alpha}{\alpha} \right)_{j} S_{r,j}^{0, 1, (\theta+k\alpha)/\alpha} \frac{(\theta +k\alpha +j\alpha +1)_{s -1}}{(\theta +k\alpha +1)_{s -1}} \frac{(\theta +k\alpha)_s (n -k\alpha)_{m-s}}{(\theta +n)_m}=
$$
$$
=\frac{(\theta +k\alpha)}{(\theta+n)_m} \sum_{j=0}^r (-1)^{r-j}\left(\frac {\theta +k\alpha +\alpha}{\alpha} \right)_{j}S_{r,j}^{0, 1, (\theta+k\alpha)/\alpha}  \sum_{s=0}^m {m \choose s}\frac{(\theta +k\alpha +j\alpha)_s (n-k\alpha)_{m-s}}{(\theta +k\alpha +j\alpha)}=
$$
$$
=\sum_{j=0}^r (-1)^{r-j}\left(\frac {\theta +k\alpha}{\alpha} \right)_{j} S_{r,j}^{0, 1, (\theta+k\alpha)/\alpha} \frac{(\theta+n +j\alpha)_m}{(\theta + n)_m}.
$$\hspace{14.7cm} $\square$\\\\
{\bf Remark 5.} As from the name {\it Beta-Binomial} distributions arise as Beta mixtures of a Binomial models, i.e. are models for the number of success in a sequence of independent trials once the probability of success has been randomized according to a Beta distribution. This, to some extent, clarifies the  proof of Proposition 1. in Favaro et al. (2009). In fact, despite they do not consider mixing explicitly over $(S_m|K_n=k)$ their proofs work in a multistep procedure that ends up in a double conditional mixing, both with a Binomial distribution and a Beta distribution.\\\\

In the next section we apply our approach to the study of the asymptotic properties of $K_m$ given $K_n$ and show how it strongly simplifies the derivation of relevant results. As a by product we obtain a  new decomposition for the limit law, different from that obtained in Favaro et al. (2009) but still a scale mixture of a Beta density and a transformation of the Mittag-Leffler density. 
For implementation of this kind of results in Bayesian nonparametrics in genomic applications, and for the need to derive asymptotic distributions connected with derivation of HPD intervals, see Favaro et al. (2009).

\subsection{Asymptotics}

We start recalling known results of which we will make use in the following. First a local limit law for the number of blocks under the $(\alpha, \theta)$ Poisson-Dirichlet model can be found e.g. in Pitman (2006). As $n \rightarrow \infty$  
\begin{equation}
\label{limitk}
\mathbb{P}_{\alpha, \theta} (K_n=k) \sim g_{\alpha, \theta}(z)n^{-\alpha}
\end{equation}
with $k \sim zn^\alpha$, where for $z >0$
\begin{equation}
\label{mittag}
g_{\alpha, \theta}(z):= \frac{\Gamma(\theta +1)}{\Gamma(\frac \theta\alpha +1)} z^{\frac \theta\alpha} g_\alpha(z),
\end{equation}
and $g_\alpha(\cdot)$ is the Mittag-Leffler density 
\begin{equation}
\label{mit}
g_\alpha(z)= {\alpha^{-1} z^{-1-1/\alpha}}{f_\alpha(z^{-1/\alpha})},
\end{equation}
for $f_\alpha(\cdot)$ the $\alpha$-stable density with $\alpha\in (0,1)$.
This implies that under $\mathbb{P}_{\alpha, \theta}$, (see Th. 8 in Pitman, 2003) almost surely and in $r$-th mean
\begin{equation}
\label{limit}
\frac{K_n}{n^{\alpha}}{\longrightarrow} Y_{\theta/\alpha}
\end{equation}
for $f_{Y_{\theta/\alpha}}(z)= g_{\alpha, \theta}(z)$. 
From  again Pitman (2006) we also know that, as $n \rightarrow \infty$, for $\alpha \in (0,1)$
\begin{equation}
\label{asinprimo}
\mathbb{E}_{\alpha, \theta}(K_n)\sim n^{\alpha} \frac{\Gamma(\theta +1)}{\alpha\Gamma(\theta +\alpha)},
\end{equation}
and for each $r > 0$ 
$$
\mathbb{E}_{\alpha, \theta} (K_n^{r}) \sim n^{\alpha r}\frac{\Gamma(\theta/\alpha +r +1) \Gamma(\theta +1)}{\Gamma(\theta +r\alpha +1) \Gamma(\theta/\alpha +1)}.
$$
It follows that for a $PD(\alpha, \theta +k\alpha)$ model we have 
\begin{equation}
\label{asinprimom}
\mathbb{E}_{\theta, \theta +k\alpha}(K_s) \sim  s^{\alpha} \frac{\Gamma(\theta +k\alpha +1)}{\alpha \Gamma(\theta +k\alpha +\alpha)},
\end{equation}
and for the $r$-th moment
\begin{equation}
\label{asinrsimo}
\mathbb{E}_{\alpha, \theta+k\alpha}(K_s^{r}) \sim s^{\alpha r } \frac{\Gamma(\theta/\alpha +k +r +1) \Gamma(\theta +k \alpha +1)}{\Gamma(\theta +k\alpha +r\alpha +1) \Gamma(\theta/\alpha +k +1)}.
\end{equation}\\
Adopting our approach to obtain a local limit for the moments of $(K_m|K_n=k)$ as in Favaro et al. (2009, Prop. 2) is just a matter to mix (\ref{asinrsimo}) over $s$ with a local limit law for $S_m|K_n=k$,
\begin{equation}
\label{localapp}
\mathbb{E}_{\alpha, \theta}(K_m^{r}|K_n=k)= \int_0^{m}\mathbb{E}_{\alpha, \theta+k\alpha} (K_s^{r}) f_{S_m|K_n=k}(s)ds.
\end{equation}
Notice that, by definition of rising factorials in terms of Gamma function,  $(x)_s=\Gamma(x+s)/\Gamma(x)$,  (\ref{polya}) may be written as 
$$
\mathbb{P}_{\alpha, \theta}(S_m =s|K_n=k) = \frac{\Gamma(\theta +n)}{\Gamma(\theta +k\alpha) \Gamma(n-k\alpha)}\frac{\Gamma(\theta +k\alpha +s)}{\Gamma(s+1)}\frac{\Gamma(n-k\alpha +m -s)}{\Gamma(m-s+1)}\frac{\Gamma(m+1)}{\Gamma(\theta +n +m)},
$$
and by Stirling approximation i.e. $\Gamma(m+a)/\Gamma(m+b) \sim m^{a-b}$ as $m \rightarrow \infty$, a local limit law for $(S_m|K=k)$, for $s \in (0,m)$, is given by
\begin{equation}
\label{localsm}
f^{\alpha, \theta}_{S_m|K_n=k}(s)\sim  \frac{\Gamma(\theta +n)}{\Gamma(\theta +k\alpha) \Gamma(n-k\alpha)} s^{\theta+k\alpha -1} (m-s)^{n-k\alpha -1} m^{-(\theta +n -1)}.
\end{equation}
\\\\
In the next Proposition we obtain the general result for $r\geq 1$ from (\ref{asinrsimo}). The case $r=1$ may be alternatively derived applying  the same operation to (\ref{asinprimom}).\\\\
{\bf Proposition 6.} Under the $(\alpha, \theta)$-Poisson-Dirichlet model the asymptotic behaviour of the $r$-th moment of $(K_m|K_n=k)$ is described by the following approximation
\begin{equation}
\label{primoasy}
\mathbb{E}_{\alpha, \theta} (K_m^{r}|K_n=k)\sim \left(\frac{\theta +k\alpha}{\alpha}\right)_r \frac{\Gamma (\theta +n)}{\Gamma(\theta +n +r\alpha)} m^{r\alpha}
\end{equation}
{\it Proof}: By (\ref{localapp}) and (\ref{localsm}) 
$$
\mathbb{E}_{\alpha, \theta} (K_m^r|K_n=k)=
$$
$$
=\frac{\Gamma(\theta/\alpha  +k +r +1)\Gamma(\theta +k\alpha +1)}{\Gamma(\theta +k\alpha +r\alpha +1) \Gamma (\theta/\alpha +k +1)} \int_0^m   \frac{\Gamma(\theta +n)}{\Gamma(\theta +k\alpha) \Gamma(n-k\alpha)} s^{\theta+k\alpha +r\alpha -1} (m-s)^{n-k\alpha -1} m^{-(\theta +n -1)}ds,
$$\\
multiplying and dividing by $\Gamma(\theta +n +r\alpha)$ it simplifies to
$$
=\left(\frac{\theta +k\alpha}{\alpha}\right)_r \frac{\Gamma(\theta +n)}{\Gamma(\theta +n +r\alpha)}  \int_0^m  \frac{\Gamma(\theta +n +r\alpha)}{\Gamma(\theta +k\alpha +r\alpha)\Gamma(n-k\alpha)} s^{\theta +k\alpha +r\alpha -1} (m-s)^{n- k\alpha -1} m^{-(\theta +n-1)} ds
$$
and by a change of variable, for $w=s/m$ and $ds=mdw$, 
$$
=\left(\frac{\theta +k\alpha}{\alpha}\right)_r \frac{\Gamma(\theta +n)}{\Gamma(\theta +n +r\alpha)}  m^{r\alpha}\int_0^1 \frac{\Gamma( \theta +n +r\alpha)}{\Gamma(\theta +k\alpha +r\alpha) \Gamma(n-k\alpha)} w^{\theta +k\alpha +r\alpha -1} (1 -w)^{n-k\alpha -1}dw=
$$
and the result follows. {\hspace{11cm}}$\square$

\bigskip
\bigskip
\noindent As for the asymptotic law of $K_m|K_n=k$, first notice that, as $m \rightarrow \infty$, we can always write
$$
\left(\frac{K_m}{m^\alpha}\Big\vert K_n=k\right){=} \left(\frac{K_m}{S_m^{\alpha}}\Big\vert K_n=k \right) *\left(\frac{S_m^\alpha}{m^\alpha}\Big\vert K_n=k\right),
$$
which can be rewritten as a product of independent random variables as
$$
\left(\frac{K_m}{m^\alpha} \Big\vert K_n=k \right){=}\left(\frac{K_m}{S_m^{\alpha}}\Big\vert  K_n=k, S_m=s\right) *\left( \frac{S_m^\alpha}{m^\alpha}\Big\vert K_n=k \right)  .
$$
Now, for the deletion of classes property of the $(\alpha, \theta)$ Poisson-Dirichlet model,
$$
\left(\frac{K_m}{S_m^{\alpha}}\Big\vert  K_n=k, S_m=s\right)_{\alpha, \theta}  \stackrel{in \hspace{0.1cm}law}{=} \left(\frac{K_s}{s^{\alpha}}\right)_{\alpha, (\theta+k\alpha)/\alpha}
$$
and by (\ref{limit}) almost surely and in $r$-th mean
\begin{equation}
\label{limks}
\frac{K_s}{s^\alpha} {\longrightarrow} Y_{(\theta +k\alpha)/\alpha}
\end{equation}
whose limit distribution, by an application of (\ref{mittag}), for $y >0$ is given by
$$
f_{Y_{(\theta +k\alpha)/\alpha}}(y)= g_{\alpha, (\theta +k\alpha)}(y)=\frac{\Gamma(\theta +k\alpha +1)}{\Gamma((\theta +k\alpha)/\alpha +1)} y^{\frac{\theta +k\alpha}{\alpha}} g_\alpha(y).
$$ \\
As for $(S_m/m|K_n=k)$, for each $m$ this is the proportion of success in a Beta-Binomial distribution of parameters $(m, \theta +k\alpha, n- k\alpha)$ to which the same asymptotic properties of the {\it P\'olya urn distribution} apply (see e.g. Johnson \& Kotz, 1977). It follows that as $m \rightarrow \infty$ almost surely
\begin{equation}
\label{limitbeta}
\left(\frac{S_m}{m}\Big\vert K_n=k\right)  {\longrightarrow} W \sim Beta(\theta +k\alpha, n -k\alpha).
\end{equation}
We are now in a position to state the following \\\\
{\bf Proposition 7.} Under the $(\alpha, \theta)$ Poisson-Dirichlet model $(K_m/m^\alpha|K_n=k)$ converges almost surely to a r.v. $Z_{n,k}^{\alpha, \theta}$ with limit distribution 
$$
f_{\alpha, \theta}^{n,k}(z)=\frac{\Gamma(\theta +n)}{\Gamma(\theta/\alpha +k) \Gamma(n -k\alpha)\alpha} z^{\theta/\alpha +k-1} 
\int_z^\infty (1- (z/v)^{1/\alpha})^{n-k\alpha -1} f_\alpha(v^{-1/\alpha})v^{-1/\alpha -1}dv,
$$
for $f_\alpha(\cdot)$ the density of the $\alpha$-stable r.v. for $\alpha \in (0,1)$.\\\\
{\it Proof}: 
By (\ref{limks}) and (\ref{limitbeta}) the almost sure limit of $(K_m/m^\alpha|K_n=k)$ exists as the product of independent r.v.s each admitting an almost sure limit, hence
$$
\left(\frac{K_m}{m^\alpha}\Big\vert K_n=k\right) \stackrel{a.s.}{\longrightarrow} Z_{n,k}^{\alpha, \theta}= Y_{(\theta +k\alpha)/\alpha} * W^\alpha.
$$
The density of  $Z_{n,k}^{\alpha, \theta}$ is given by
$$
f_Z (z)=\int_0^1 f_Y(z w^{-\alpha})w^{-\alpha} f_W(w) dw=
$$
$$
=\frac{\Gamma(\theta +k\alpha +1)}{\alpha\Gamma(\theta/\alpha+k+1)}\int_0^1 (zw^{-\alpha})^{\theta/\alpha +k -1 -1/\alpha} f_\alpha[(zw^{-\alpha})^{-1/\alpha}] \frac{1}{w^{\alpha}} \frac{\Gamma(\theta +n)}{\Gamma(\theta +k\alpha)\Gamma(n- k\alpha)} w^{\theta +k\alpha -1}(1 -w)^{n-k\alpha -1}dw=
$$\\
which simplifies to
$$
=\frac{\Gamma(\theta +n)}{\Gamma(\theta/\alpha +k)\Gamma(n-k\alpha)} z^{\theta/\alpha +k -1-1/\alpha} \int_0^1 (1 -w)^{n- k\alpha -1} f_\alpha [(zw^{-\alpha})^{-1/\alpha}]dw=
$$
and by the change of variable $zw^{-\alpha}=v$, $w=(zv^{-1})^{1/\alpha}$, $dw=\alpha^{-1} z^{1/\alpha}v^{-1/\alpha -1}dv$, it follows
$$
=\frac{\Gamma(\theta +n)}{\Gamma(\theta/\alpha +k)\Gamma(n-k\alpha)\alpha} z^{\theta/\alpha +k -1} \int_0^\infty (1 -(zv^{-1})^{1/\alpha})^{n -k\alpha -1} f_\alpha (v^{-1/\alpha}) v^{-1/\alpha -1}dv.
$$\hspace{ 15cm}$\square$\\\\
Next Proposition proves both the convergence in $r$-th mean of $(K_m/m^\alpha|K_n=k)$ to $Z_{n,k}^{\alpha, \theta}$ and that our result, while agrees with Favaro et al. (2009, Proposition 2.) provides a new decomposition for the limit law.\\\\
{\bf Proposition 8.} Let $H=Y_1*X$  for  $Y_1$ and $X$ independent r.v.s, $Y_1 \sim g_{\alpha, (\theta+n)}$ and $X\sim Beta(\theta/\alpha +k,  n/\alpha -k)$, then  $Z_{n,k}^{\alpha, \theta}$ and $H$ have the same characteristic function
$$
G_{\alpha, \theta}^{n,k}(t)= \sum_{r \geq 0} \frac{(it)^r}{r!} \left(\frac{\theta +k\alpha}{\alpha}\right)_r
 \frac{1}{(\theta +n)_{r\alpha}}.
$$
\\
{\it Proof}:   First notice that Proposition 3. is enough to say that for $m \rightarrow \infty$
$$
\mathbb{E}_{\alpha, \theta }\left(\frac{K_m^r}{m^{r\alpha}}\Big\vert  K_n=k\right) \stackrel{}{\longrightarrow} \left(\frac{\theta +k\alpha}{\alpha}\right)_r \frac{\Gamma (\theta +n)}{\Gamma(\theta +n +r\alpha)}.
$$
Now the density of  $Z_{n,k}^{\alpha, \theta}$ may be written as
$$
f_Z(z)=\frac{\Gamma(\theta +n)\Gamma(\theta+k\alpha+1)}{\Gamma(\theta +k\alpha)\Gamma(n -k\alpha) \Gamma((\theta  +k\alpha)/\alpha +1)}\frac 1\alpha z^{\theta/\alpha+k-1}\int_z^{\infty}\alpha^{-1}s^{-1/\alpha-1} f_\alpha(s^{-1/\alpha})\left( 1-(z/s)^{1/\alpha}\right)^{n-k\alpha-1}ds
$$
whose characteristic function by (\ref{mittag}) and (\ref{mit}) is given by 
$$
=\frac{\Gamma(\theta +n)\Gamma(\theta+k\alpha+1)}{\Gamma(\theta +k\alpha)\Gamma(n -k\alpha) \Gamma((\theta  +k\alpha)/\alpha+1)} \frac 1\alpha \int_0^\infty \exp\{itz\} z^{\theta/\alpha+k-1}\int_z^{\infty}g_{\alpha}(s)\left( 1-(z/s)^{1/\alpha}\right)^{n-k\alpha-1}ds dz
$$
this may be rewritten as
$$
=\frac{\Gamma(\theta+k\alpha+1)}{\Gamma((\theta  +k\alpha)/ \alpha+1)} 
\frac 1\alpha \int_z^{\infty}g_{\alpha}(s)
\int_0^\infty \exp\{itz\}\frac{ \Gamma(\theta +n)}{ \Gamma(\theta +k\alpha)\Gamma(n -k\alpha)}z^{\theta/\alpha+k-1}
\left( 1-(z/s)^{1/\alpha}\right)^{n-k\alpha-1}dz ds \nonumber
$$
and by a change of variable $(z/s)^{1/\alpha}=y$, $z=y^{\alpha}s$, $dz=s\alpha y^{\alpha-1}dy$ reduces to
$$
=\frac{\Gamma(\theta+k\alpha+1)}{ \Gamma((\theta  +k\alpha)/\alpha+1)} 
\frac 1\alpha \int_0^{\infty} g_{\alpha}(s)
\int_0^s e^{ity^{\alpha}s}\frac{\Gamma(\theta +n)}{ \Gamma(\theta +k\alpha)\Gamma(n -k\alpha)}(y^{\alpha}s)^{\theta/\alpha+k-1}
\left( 1- y\right)^{n-k\alpha-1}s \alpha y^{\alpha -1}dy ds
$$
and then to
$$
=\frac{\Gamma(\theta+k\alpha+1)}{ \Gamma((\theta  +k\alpha)/\alpha +1)} 
\int_0^{\infty} s^{\theta/\alpha +k} g_{\alpha}(s)
\int_0^1 e^{ity^{\alpha}s}\frac{\Gamma(\theta +n)}{ \Gamma(\theta +k\alpha)\Gamma(n -k\alpha)}(y)^{\theta+k\alpha-1}
\left( 1- y\right)^{n-k\alpha-1}dy ds.
$$
By the characteristic function of $Y^\alpha$ for $Y \sim Beta(\theta +k\alpha, n -k\alpha)$ we can write
$$
\label{last}
=\frac{\Gamma(\theta+k\alpha+1)}{ \Gamma((\theta  +k\alpha)/\alpha +1)} \sum_{r=0}^{\infty} \frac{(it)^{r}}{r!}\frac{(\theta +k\alpha)_{r\alpha} }{(\theta +n)_{r\alpha}}  \int_0^{\infty} s^{\theta/\alpha +k +r} g_{\alpha}(s) ds
$$
and by (\ref{mittag})
\begin{equation}
\label{last}
=\sum_{r=0}^{\infty} \frac{(it)^{r}}{r!}\frac{(\theta +k\alpha)_{r\alpha }}{(\theta +n)_{r\alpha}} \frac{\Gamma(\theta+k\alpha+1)}{ \Gamma((\theta  +k\alpha)/\alpha +1)}  \frac{\Gamma((\theta +k\alpha +r\alpha)/\alpha +1)}{\Gamma(\theta +k\alpha +r\alpha +1)}.
\end{equation}
By the usual properties of Gamma function the last expression corresponds to
$$
=\sum_{r=0}^{\infty} \frac{(it)^{r}}{r!}\frac{\Gamma(\theta + k\alpha +r\alpha)\Gamma(\theta +n)}{\Gamma(\theta +k\alpha) \Gamma(\theta +n +r\alpha)}\frac{(\theta +k\alpha)\Gamma(\theta +k\alpha)}{\frac{\theta +k\alpha}{\alpha}\Gamma(\frac{\theta +k\alpha}{\alpha}) }\frac{\Gamma(\frac{\theta +k\alpha}{\alpha}+r) \frac{\theta +k\alpha +r\alpha}{\alpha}}{(\theta +k\alpha +r\alpha) \Gamma(\theta +k\alpha +r\alpha)}
$$
which simplifies to
$$
=\sum_{r=0}^{\infty} \frac{(it)^{r}}{r!}\left(\frac{\theta +k\alpha}{\alpha}\right)_r\frac{1}{(\theta +n)_{r\alpha}}
$$
and the conclusion follows by the result in Proposition 2 in Favaro et al. (2009).

\section*{Acknowledgement} The author wishes to thank Andrea Tancredi and Serena Arima for some helpful and illuminating discussions.

\section*{References}
\newcommand{\bibu}{\item \hskip-1.0cm}
\begin{list}{\ }{\setlength\leftmargin{1.0cm}}

\bibu \textsc {Cerquetti, A.} (2008) Generalized Chinese restaurant construction of exchangeable Gibbs partitions and related results. {\sf http://arxiv.org/abs/0805.3853}

\bibu  \textsc {Charalambides, C. A.} (2005) {\it Combinatorial Methods in Discrete Distributions}. Wiley, Hoboken NJ.

\bibu \textsc{Favaro, S., Lijoi, A., Mena, R. and Pr\"unster, I.} (2009) Bayesian non-parametric inference for species variety with a two-parameter Poisson-Dirichlet process prior. {\it JRSS-B}, 71, 993-1008.

\bibu \textsc{Gnedin, A. and Pitman, J. } (2006) {Exchangeable Gibbs partitions  and Stirling triangles.} {\it Journal of Mathematical Sciences}, 138, 3, 5674--5685. 

\bibu \textsc{Gnedin, A., Haulk, S. and Pitman, J.} (2009) {Characterizations of exchangeable partitions and random discrete distributions by deletion properties}. {\sf http://arxiv.org/abs/0909.3642}



\bibu \textsc{Hsu, L. C, \& Shiue, P. J.} (1998) A unified approach to generalized Stirling numbers. {\it Adv. Appl. Math.}, 20, 366-384.

\bibu \textsc{Johnson, N.L. \& Kotz, S.} (1977) {\it Urn models and their application}. Wiley \& Sons.

\bibu \textsc{Johnson, N. L. \& Kotz, S.} (2005) {\it Univariate Discrete Distributions}. 3rd Edition, Wiley \& Sons.

\bibu \textsc{Lijoi, A., Mena, R. and Pr\"unster, I.} (2007) Bayesian nonparametric estimation of the probability of discovering new species.  {\it Biometrika}, 94, 769--786.

\bibu \textsc{Lijoi, A., Pr\"unster, I. and Walker, S.G.} (2008) Bayesian nonparametric estimator derived from conditional Gibbs structures. {\it Annals of Applied Probability}, 18, 1519--1547.


\bibu \textsc{Perman, M., Pitman, J, \& Yor, M.} (1992) Size-biased sampling of Poisson point processes and excursions. {\it Probab. Th. Rel. Fields}, 92, 21--39.

\bibu \textsc{Pitman, J.} (1995) Exchangeable and partially exchangeable random partitions. {\it Probab. Th. Rel. Fields}, 102: 145-158.

\bibu \textsc{Pitman, J.} (1996a) Some developments of the Blackwell-MacQueen urn scheme. In T.S. Ferguson, Shapley L.S., and MacQueen J.B., editors, {\it Statistics, Probability and Game Theory}, volume 30 of {\it IMS Lecture Notes-Monograph Series}, pages 245--267. Institute of Mathematical Statistics, Hayward, CA.

\bibu \textsc{Pitman, J.} (1996b) Notes on the two parameter generalization of Ewens random partition structure. {\it Manuscript} University of California, Berkeley. Unpublished.

\bibu \textsc{Pitman, J.} (2003) {Poisson-Kingman partitions}. In D.R. Goldstein, editor, {\it Science and Statistics: A Festschrift for Terry Speed}, volume 40 of Lecture Notes-Monograph Series, pages 1--34. Institute of Mathematical Statistics, Hayward, California.

\bibu \textsc{Pitman, J.} (2006) {\it Combinatorial Stochastic Processes}. Ecole d'Et\'e de Probabilit\'e de Saint-Flour XXXII - 2002. Lecture Notes in Mathematics N. 1875, Springer.

\bibu \textsc{Pitman, J. and Yor, M.} (1997) The two-parameter Poisson-Dirichlet distribution derived from a stable subordinator. {\it Ann. Probab.}, 25:855--900.


\bibu \textsc{Yamato, H. and Sibuya, M.} (2000) Moments of some statistics of Pitman sampling formula. {\it Bull. Inform. Cybernet.}, 32 1--10.

\end{list}
\end{document}